# EFFECTIVE REDSHIFT

## TRISTAN YANG

Abstract. The "higher chromatic" Quillen-Lichtenbaum conjecture, as proposed by Ausoni and Rognes, posits that the finite localization map $K(R) \to L_{n+1}^f K(R)$ is a $p$-local equivalence in large degrees for suitable ring spectra $R$. We give a simple criterion in terms of syntomic cohomology for an effective version of Quillen-Lichtenbaum, i.e. for identifying the degrees in which the localization map is an isomorphism. Combining our result with recent computations implies that the finite localization map is $(-1)$-truncated in the cases $R = \mathrm{BP}\langle n \rangle$, $R = k(n)$, and $R = \mathrm{ko}$.

The classical *Quillen-Lichtenbaum conjecture* [Qui75], now proved by Voevodsky, Rost, and others [Voe03; Voe11], posits that algebraic $K$-theory should be computable in terms of étale $K$-theory. This admits a reinterpretation in the language of chromatic homotopy theory: by a result of Thomason [Tho82], the $p$-completed localization $L_1 K(R)_p^\wedge$ of $K(R)$ at complex $K$-theory satisfies étale descent, and the Quillen-Lichtenbaum conjecture is equivalent to asking that the map

$$K(R) \to L_1 K(R)$$

be a $p$-local equivalence in large degrees for suitable rings $R$ (see [Wal84, §4]).

From this perspective, one is invited to consider a "higher" analogue by taking $R$ to be a suitable homotopy ring spectrum and replacing $L_1$ with one of successive higher chromatic localizations $L_n^f$. This was first formulated by Ausoni and Rognes [AR02]:

**Hope 0.1.** For suitable rings $R$ of "chromatic height" $n$, the map

$$K(R) \to L_{n+1}^f K(R)$$

is a $p$-local equivalence in large degrees.

Hope 0.1 realizes the philosophy that $K$-theory increases chromatic height by +1, giving rise to the name *"chromatic redshift"*. Beginning with [AR02], who computed redshift for $R = \ell_p^\wedge$, various redshift results have appeared in [Aus10; AR12; AACHR22; LMMT22; CMNN22]. Recently, Hahn and Wilson [HW22] demonstrated the first arbitrary-height examples of redshift for the case that $R$ is a truncated Brown-Peterson spectra $\mathrm{BP}\langle n \rangle$, which are each of height $n$ in the sense of the *fp-type* of [MR99]. Additionally, [BSY22] proved a different weak form of redshift for all $\mathbb{E}_\infty$-rings $R$, showing that if $L_{T(n)} R \neq 0$ then $L_{T(n+1)} K(R) \neq 0$.

In this work, we aim to make Hope 0.1 *effective*; that is, find the range in which the map $K(R) \to L_{n+1}^f K(R)$ is an isomorphism. We first show in Section 1 that it suffices to bound the homotopy groups of $\mathrm{TC}(R) \otimes F$, where $F$ is a type-$(n + 2)$ finite complex, in terms of the dimension of $F$:

**Theorem 1.3.** *Let $X$ be a $p$-local spectrum, and suppose there exists a type-$(n+2)$ finite spectrum $F$ such that $X \otimes F$ is $(\dim F + d)$-truncated. Then the localization map*

$$X \to L_{n+1}^f X$$

*is d-truncated (i.e. the fiber is concentrated in degrees $\leq d$).*

In [Section 2](#) we show how the machinery of syntomic cohomology and the motivic spectral sequence, together with the results of [Section 1](#), can be used to prove effective redshift in several cases. In particular, recent computations of the syntomic cohomology of ko [AAR23], BP⟨n⟩ [AHW25], and $k(n)$ [AHW24] imply:

**Example 2.9.** *The localization maps*

$$K(\text{ko})_{(2)} \to L_{n+1}^f K(\text{ko})_{(2)}$$
$$K(\text{ko})_{(2)} \to L_{n+1} K(\text{ko})_{(2)}$$

*are $(-1)$-truncated. In particular, they are an isomorphism on $\pi_*$ for $* > 0$.*

**Theorem 2.7.** *Suppose we are given a form of* BP⟨n⟩ *whose syntomic cohomology is given by [Example 2.6](#), as in [HRW25; AHW25]. Then the localization maps*

$$K(\text{BP}\langle n \rangle)_{(p)} \to L_{n+1}^f K(\text{BP}\langle n \rangle)_{(p)}$$
$$K(\text{BP}\langle n \rangle)_{(p)} \to L_{n+1} K(\text{BP}\langle n \rangle)_{(p)}$$

*are $(-1)$-truncated. In particular, they are an isomorphism on $\pi_*$ for $* > 0$.*

**Theorem 3.13.** *For any $\mathbb{E}_1$ form of $k(n)$, the localization maps*

$$K(k(n))_{(p)} \to L_{n+1}^f K(k(n))_{(p)}$$
$$K(k(n))_{(p)} \to L_{n+1} K(k(n))_{(p)}$$

*are $(-1)$-truncated.*

Therefore, the Quillen-Lichtenbaum range for BP⟨n⟩ and for $k(n)$ are in fact "as good as possible."

Towards exploring more general conditions under which redshift holds, we provide in [Section 3](#) certain criteria for the motivic associated graded of THH($R$) that reduce redshift to the *Segal Conjecture* ([Theorem 3.6](#)). As a rudimentary application, we then show how this can help simplify somewhat the proof of effective redshift for $k(n)$.

## 1. EFFECTIVE QUILLEN-LICHTENBAUM

In [HW22], the starting point for proving Quillen-Lichtenbaum-type statements is the following result of Mahowald and Rezk [MR99, Theorem 8.2(2)].

**Proposition 1.1.** *Let $X$ be a $p$-local spectrum, and suppose there exists a type-$(n + 1)$ finite spectrum $F$ such that $X \otimes F$ is truncated. Then the localization map*

$$X \to L_n^f X$$

*is truncated.*

**Remark 1.2.** In [MR99], attention is restricted to $fp$-spectra $X$; that is, $p$-complete bounded below spectra whose $\mathbb{F}_p$-cohomology is finitely presented over the mod $p$ Steenrod algebra. We note here, though, that the result as stated above holds in generality for all $p$-local spectra $X$, and we give a simple self-contained proof below. This also appears as [BHLS23, Lemma 7.22].

*Proof of Proposition 1.1.* The collection of all $F$ such that $X \otimes F$ is truncated form a thick subcategory. Therefore, by the thick subcategory theorem, it suffices to consider the case that $F$ is a generalized Moore spectrum $M(n) \coloneqq \mathbb{S}/(p^{i_0}, v_1^{i_1}, \ldots, v_n^{i_n})$. Furthermore, note that $\mathrm{fib}(X \to L_n^f X) \otimes F$ is also truncated, as it is in fact equivalent to $X \otimes F$. In other words, the fiber of the localization maps also satisfies the hypotheses of the proposition, so by substituting the fiber in the place of $X$ we may assume without loss of generality that $L_n^f X \simeq 0$.

We now claim that, in fact, $X$ is still truncated after tensoring only with the type-$n$ finite spectrum $M(n-1)$. Indeed, consider the cofiber sequence

$$\Sigma^{i_n |v_n|} M(n-1) \otimes X \xrightarrow{v_n^{i_n}} M(n-1) \otimes X \longrightarrow M(n) \otimes X.$$

By the induced long exact sequence in homotopy, the self-map $v_n^{i_n}$ is an isomorphism on $\pi_k(M(n-1) \otimes X)$ for large degrees $k$. Therefore, for large $k$ we have

$$\pi_k(M(n-1) \otimes X) \cong \pi_k(M(n-1)[v_n^{-1}] \otimes X)$$

but the right-hand side is exactly a version of $\pi_k(T(n) \otimes X)$, which is equivalent to $\pi_k(T(n) \otimes L_{T(n)} X) = 0$.

After iterating the above argument, we eventually deduce that $X$ is truncated after tensoring with the type-0 spectrum $\mathbb{S}_{(p)}$; i.e. $X$ itself is truncated. □

We now show that the above proposition can be strengthened to an effective version. For $F$ a $p$-local finite spectrum, let $\dim(F)$ refer to the highest nonzero degree of $H_*(F; \mathbb{F}_p)$.

**Theorem 1.3.** *Let $X$ be a $p$-local spectrum, and suppose there exists a type-$(n+1)$ finite spectrum $F$ such that $X \otimes F$ is $(\dim F + d)$-truncated. Then the localization map*

$$X \to L_n^f X$$

*is $d$-truncated.*

*Proof.* Again, by replacing $X$ with the fiber of the localization map, we assume without loss of generality that $L_n^f X \simeq 0$. By Proposition 1.1, $X$ is truncated. Let $d_X$ be the degree of its top homotopy group, so we wish to show that $d_X \leq d$. We proceed via the Atiyah-Hirzebruch spectral sequence

$$H_s(F; \pi_t X) \implies \pi_{s+t}(X \otimes F)$$

For each $t$, we claim that $H_*(F; \pi_t X)$ has top degree $* = \dim F$. To show this, note that $L_{T(0)} X \simeq 0$, so $\pi_* X$ is $p$-power torsion. By a theorem of Kulikov [Rot95, Theorem 10.36], such a group is an extension of a direct sum of copies of $\mathbb{Z}/p^k$ by a direct sum of copies of $\mathbb{Z}/p^\infty$, so it suffices to consider the top degree of $H_*(F; \mathbb{Z}/p^k)$ for each $k \leq \infty$. For finite $k$, this is straightforward from the iterated cofiber sequences

$$H\mathbb{Z}/p^k \otimes F \xrightarrow{p} H\mathbb{Z}/p^{k+1} \otimes F \longrightarrow H\mathbb{F}_p \otimes F.$$

But observe that the maps $H_*(F; \mathbb{Z}/p^k) \to H_*(F; \mathbb{Z}/p^{k+1})$ must be injective in the top degree, and therefore colimit $H_{\dim F}(F; \mathbb{Z}/p^\infty)$ also does not vanish.

Now in the Atiyah-Hirzebruch spectral sequence, the group in location $(s, t) = (\dim F, d_X)$ on the $E_2$-page is not affected by any differentials. Therefore, the top homotopy group of $X \otimes F$ is in degree $(\dim F + d_X)$, so $d_X \leq d$ as desired. □

We can now already apply Theorem 1.3 to redshift computations that have been performed at low heights.

**Example 1.4.** Let $X = K(R)_{(p)}$ for $R = \mathrm{ku}_p, \mathrm{ku}_{(p)}, \ell_p, \ell$, or $k(1)$. Let $p \geq 5$, and let $F = V(1)/v_2^k$ where $V(1)$ denotes the Smith-Toda complex $\mathbb{S}/(p, v_1)$. By the results of [AR12, Theorem 1.1], [Aus10, Theorem 8.1], [BM08], and [AR02, Theorem 9.1], $F \otimes X$ vanishes in degrees

$$* > k(2p^2 - 2) + (2p - 2) + 2 = \dim F - 1.$$

Therefore, by Theorem 1.3, the localization map

$$K(R)_{(p)} \to L_2^f K(R)_{(p)}$$

is $(-1)$-truncated.

**Example 1.5.** Let $X = K(\mathrm{BP}\langle 2 \rangle)_{(p)}$, and let $F = V(2)/v_3^k$ where $p \geq 7$. By [AACHR22, Theorem 1.1], $F \otimes X$ vanishes in degrees

$$* > k(2p^3 - 2) + (2p^2 - 2) + (2p - 2) + 3 = \dim F - 1.$$

Therefore, by Theorem 1.3, the localization map

$$K(\mathrm{BP}\langle 2 \rangle)_{(p)} \to L_3^f K(\mathrm{BP}\langle 2 \rangle)_{(p)}$$

is $(-1)$-truncated.

To replicate Example 1.4 and Example 1.5 at higher heights, we run into two problems. First, it becomes progressively harder to calculate $V(n)_* \mathrm{TC}(R)$ directly, as done for the results cited above. Secondly, existence of the Smith-Toda complex $V(n)$ requires the assumption that $p$ be increasingly large as $n$ grows. We address both of these issues in the next section. We will also be able to see that the conclusions of Example 1.4, for $R = \ell = \mathrm{BP}\langle 1 \rangle$, and Example 1.5 are still true without the restriction on $p$.

## 2. Redshift and the Motivic Filtration

The goal of this section is to show how Theorem 1.3 can be combined with recent work of Hahn, Raksit, and Wilson [HRW25] to yield wider effective redshift results for TC and $K$-theory. As suggested by Example 1.4 and Example 1.5, we would ideally like to have a Smith-Toda complex in order to apply Theorem 1.3 most easily—but unfortunately, Smith-Toda complexes do not exist in general. The key idea is to instead mimic the approach after passing to the world of $\mathbb{C}$-motivic spectra.

Specifically, given an algebra $R$ over a suitable $\mathbb{E}_\infty$-ring $A$, Hahn, Raksit and Wilson [HRW25] construct a motivic filtration $\mathrm{fil}^*_{\mathrm{mot}/A} \mathrm{TC}(R)_p^\wedge$ on the $p$-completed topological cyclic homology of $R$. By work of Gheorghe-Isaksen-Krause-Ricka [GIKR22], this can be viewed as a $\mathbb{C}$-motivic spectrum that recovers $\mathrm{TC}(R)_p^\wedge$ upon inverting $\tau$, while the associated graded $\mathrm{gr}^*_{\mathrm{mot}/A} \mathrm{TC}(R)_p^\wedge$ (syntomic cohomology) is a module over $\mathrm{gr}^*_{\mathrm{ev}} \mathbb{S}$, the cofiber of $\tau$. Now, the category of $\mathrm{gr}^*_{\mathrm{ev}} \mathbb{S}$-modules admits a $t$-structure whose heart is equivalent to the category of even $\mathrm{MU}_* \mathrm{MU}$-comodules [Pst23; GIKR22]. In particular, there is a $\mathrm{gr}^*_{\mathrm{ev}} \mathbb{S}$-module $\mathrm{gr}^*_{\mathrm{ev}} \mathbb{S}/(p, \ldots, v_k)$ corresponding to the comodule $\mathrm{MU}_*/(p, \ldots, v_k)$ under this equivalence. Following [HRW25], we use the general notation

$$(-)/(p, \ldots, v_k) := (-) \underset{\mathrm{gr}^*_{\mathrm{ev}} \mathbb{S}}{\otimes} \mathrm{gr}^*_{\mathrm{ev}} \mathbb{S}/(p, \ldots, v_k).$$

This can be thought of as an algebraic analogue of tensoring with a Smith-Toda.

We now port Theorem 1.3 to this setting. Let us abbreviate

$$m_{n+1} := \sum_{i=0}^{n+1} (|v_i| + 1)$$

the would-be dimension of a hypothetical Smith-Toda $V(n + 1)$.

**Remark 2.1.** We note that there is possible linguistic ambiguity stemming from the fact that the homotopy groups of $\mathrm{gr}^\star_{\mathrm{ev}}\mathbb{S}$-modules are doubly graded: by "$d$-truncated" we will mean to refer to the topological degree—i.e. that $\pi_*\mathrm{gr}^k_{\mathrm{ev}}$ vanishes in degrees $* > d$ for each graded piece $k$.

**Remark 2.2.** Note that after modding out by a power of $p$, we are free to drop the $p$-completion decoration unambiguously from our notation for TC.

**Theorem 2.3.** *Let $R$ be an algebra over a chromatically $p$-quasisyntomic ring $A$, in the sense of [HRW25, Definition 1.3.1], and suppose that the mod-$(p,\dots,v_{n+1})$ syntomic cohomology $(\mathrm{gr}^\star_{\mathrm{mot}/A}\mathrm{TC}(R))/(p,\dots,v_{n+1})$ of $R$ is $(m_{n+1} + d)$-truncated. Then the localization map*

$$\mathrm{TC}(R) \to L^f_{n+1}\mathrm{TC}(R)$$

*is $d$-truncated.*

*Proof.* Let $M(n+1) = \mathbb{S}/(p^{i_0},\dots,v_{n+1}^{i_{n+1}})$ denote a generalized Moore complex, so that we have a motivic spectral sequence

$$\pi_*(\mathrm{gr}^\star_{\mathrm{mot}/A}\mathrm{TC}(R)/(p^{i_0},\dots,v_{n+1}^{i_{n+1}})) \implies M(n+1)_*\mathrm{TC}(R).$$

We can resolve the above $E_2$-page by finitely many copies of $\mathrm{gr}^\star_{\mathrm{mot}/A}\mathrm{TC}(R)/(p,\dots,v_{n+1})$ via the cofiber sequences

$$\mathrm{gr}^\star_{\mathrm{mot}/A}\mathrm{TC}(R)/(\dots,v_m^k,\dots) \xrightarrow{v_m} \mathrm{gr}^\star_{\mathrm{mot}/A}\mathrm{TC}(R)/(\dots,v_m^{k+1},\dots)$$
$$\downarrow$$
$$\mathrm{gr}^\star_{\mathrm{mot}/A}\mathrm{TC}(R)/(p,\dots,v_{n+1}).$$

If the leftmost term vanishes in degrees $> d_0$, then the middle term vanishes in degrees $> d_0 + |v_m|$. Therefore, we get that $\mathrm{gr}^\star_{\mathrm{mot}/A}\mathrm{TC}(R)/(p,\dots,v_{n+1}^{i_{n+1}})$ vanishes in degrees greater than

$$d + m_{n+1} + \sum_{k=0}^{n+1}(i_k - 1)(2p^k - 2) = d + (n+2) + \sum_{k=0}^{n+1} i_k(2p^k - 2)$$
$$= d + \dim M(n+1).$$

The vanishing bound on the $E_2$ page of the aforementioned spectral sequence implies the same bound on $M(n+1)_*\mathrm{TC}(R)$, so the result follows from Theorem 1.3. $\qquad\square$

As an immediate corollary of the Dundas-Goodwillie-McCarthy Theorem [DGM12], we get a result for $K$-theory.

**Corollary 2.4.** *Let $R$ be an algebra over a chromatically $p$-quasisyntomic $\mathbb{E}_\infty$-ring $A$, and suppose that $\mathrm{TC}(\pi_0 R)^\wedge_p$ and the mod-$(p,\dots,v_{n+1})$ syntomic cohomology of $R$ are $(m_{n+1} + d)$-truncated. If the localization map*

$$K(\pi_0 R)_{(p)} \to L^f_{n+1}K(\pi_0 R)_{(p)}$$

*is $d$-truncated, then the localization map*

$$K(R)_{(p)} \to L^f_{n+1}K(R)_{(p)}$$

*is also $d$-truncated.*

*Proof.* This follows from Theorem 2.3, the ($p$-completed) Dundas-Goodwillie-McCarthy pullback square

$$\begin{CD}
K(R)^\wedge_p @>>> \mathrm{TC}(R)^\wedge_p \\
@VVV @VVV \\
K(\pi_0 R)^\wedge_p @>>> \mathrm{TC}(\pi_0 R)^\wedge_p
\end{CD}$$

and the arithmetic pullback square

$$\begin{CD}
X_{(p)} @>>> X^\wedge_p \\
@VVV @VVV \\
X_{(p)}[p^{-1}] @>>> X^\wedge_p[p^{-1}]
\end{CD}$$

where we note that the localization map is an equivalence on rational spectra. $\square$

Let us also remark on a criterion that allows one to replace the finite $L^f_n$-localization with $L_n$-localization. Recall that the Adams grading on the motivic spectral sequence places the group $\pi_t \mathrm{gr}^s(-)$ in column and row $(t, 2s - t)$. Whereas the hypothesis of Theorem 2.3 asks for a vertical vanishing line, we can also ask for a horizontal one:

**Proposition 2.5.** *Let $R$ be an algebra over a chromatically $p$-quasisyntomic $\mathbb{E}_\infty$-ring $A$, and suppose the motivic spectral sequence for $F_* \mathrm{TC}(R)$ has an eventual horizontal vanishing line when displayed in Adams grading, where $F$ is some type-$(n + 2)$ finite complex. Then the height-$(n + 1)$ telescope conjecture is true of $\mathrm{TC}(R)$, that is, the map*

$$L^f_{n+1}\mathrm{TC}(R) \to L_{n+1}\mathrm{TC}(R)$$

*is an equivalence.*

*Proof.* We duplicate the argument of [HRW25, Theorem 6.6.4]. It suffices to show the result for $\mathrm{TC}^-(R)$ and $\mathrm{TP}(R)$; the argument is the same for both so we record only the case $\mathrm{TC}^-(R)$ below. (Alternatively, note that the latter is a module over the former and the localizations are smashing.)

First, for each $m \le n + 1$, to check that

$$L_{T(m)}\mathrm{TC}^-(R) \to L_{K(m)}\mathrm{TC}^-(R)$$

is an equivalence we check that it is a $T(m)$-equivalence (note that on the left hand side we have $T(m) \otimes L_{T(m)}\mathrm{TC}^-(R) \simeq T(m) \otimes \mathrm{TC}^-(R)$). Let us set $T(m) = M(m - 1)[v_m^{-1}]$ for $M(m - 1) = \mathbb{S}/(p^{i_0}, \ldots, v_{m-1}^{i_{m-1}})$ a generalized Moore complex. By assumption, the motivic spectral sequence for $M(n + 1)_* \mathrm{TC}(R)$ has a horizontal vanishing line. But note that the $v_{n+1}$-Bockstein spectral sequence computing the motivic $E_2$-page for $M(n)_* \mathrm{TC}(R)$ has an $E_2$ page obtained by adjoining an element $v_{n+1}$ in Adams weight 0. Therefore, the motivic spectral sequence for $M(n)_* \mathrm{TC}(R)$ also has a horizontal vanishing line, and similarly for lower heights.

To complete the proof, recall that the motivic filtration can be realized as

$$\mathrm{fil}^\bullet_{\mathrm{mot}/A} \mathrm{TC}^-(R) \simeq \lim_\Delta \tau_{\ge 2*} C^\bullet$$

where $C^\bullet$ is the complex associated to descent along some appropriate map $\mathrm{THH}(A) \to \mathrm{THH}(A \otimes \mathrm{MU}/\mathbb{S})$, as described in [HW22, Remark 1.3.8]. There is an equivalence

$$M(m - 1) \otimes \mathrm{TC}^-(R) \simeq \lim_\Delta (M(m - 1) \otimes C^\bullet)$$

where we use the finiteness of $M(m-1)$ to pass it inside the totalization. Now, because the motivic spectral sequence coincides with the (double-speed décalage of) the descent spectral sequence, by [CM21, Lemma 2.34] and [Mat15, Proposition 3.12] we may also pass the $v_m^{-1}$ inside the totalization, to obtain

$$T(m) \otimes \mathrm{TC}^-(R) \simeq \lim_\Delta (T(m) \otimes C^\bullet).$$

We thus reduce to checking the telescope conjecture on each term in the totalization, for which it holds because each $C^i$ is an MU-module. $\qquad\square$

**Example 2.6.** Suppose we have a form of $\mathrm{BP}\langle n \rangle$ such that

$$\mathrm{gr}^\star_{\mathrm{mot/MU}} \mathrm{TC}(\mathrm{BP}\langle n \rangle)/(p,\dots,v_{n+1})$$

$$\cong \Lambda_{\overline{\mathbb{F}}_p}(\partial, \lambda_1, \dots, \lambda_{n+1}) \oplus \bigoplus_{k=1}^{n+1} \Lambda_{\overline{\mathbb{F}}_p}(\lambda_1, \dots, \hat{\lambda}_k, \dots, \lambda_{n+1}) \otimes \mathbb{F}_p\{\lambda_k t^{kd} : 0 < d < p\}$$

where $|\partial| = (-1,1)$, $|\lambda_k| = (2p^k - 1, 1)$ and $|t| = (-2, 0)$. This computation at height $n = 1$ was recently performed in [HRW25, Theorem 6.0.2], and the case for higher heights is to appear in forthcoming work [AHW25].

By Theorem 2.3 and Proposition 2.5, we obtain in this situation that the localization maps

$$\mathrm{TC}(\mathrm{BP}\langle n \rangle) \to L_{n+1}^f \mathrm{TC}(\mathrm{BP}\langle n \rangle)$$

$$\mathrm{TC}(\mathrm{BP}\langle n \rangle) \to L_{n+1} \mathrm{TC}(\mathrm{BP}\langle n \rangle)$$

are $(-1)$-truncated. In particular, they are an isomorphism on $\pi_*$ in positive degrees.

We deduce:

**Theorem 2.7.** *In the situation of Example 2.6, the localization maps*

$$K\langle \mathrm{BP}\langle n \rangle)_{(p)} \to L_{n+1}^f K(\mathrm{BP}\langle n \rangle)_{(p)}$$

$$K(\mathrm{BP}\langle n \rangle)_{(p)} \to L_{n+1} K(\mathrm{BP}\langle n \rangle)_{(p)}$$

*are $(-1)$-truncated. In particular, they are an isomorphism on $\pi_*$ in positive degrees.*

*Proof.* Corollary 2.4 reduces us to showing the statement for $\pi_0 \mathrm{BP}\langle n \rangle \cong \mathbb{Z}_{(p)}$, for which we recapitulate classical Quillen-Lichtenbaum results. By [CM21, Theorem 6.18], the map

$$K(\mathbb{Z}_{(p)})_{(p)} \to K^{\mathrm{Sel}}(\mathbb{Z}_{(p)})_{(p)}$$

is an isomorphism in degrees $\geq \mathrm{vcd}_p(\mathbb{Q}) - 2 = 0$, where Selmer $K$-theory $K^{\mathrm{Sel}}$ is defined via the pullback

$$\begin{array}{ccc} K^{\mathrm{Sel}}(R)_{(p)} & \longrightarrow & \mathrm{TC}(R)_{(p)} \\ \downarrow & & \downarrow \\ L_1 K(R)_{(p)} & \longrightarrow & L_1 \mathrm{TC}(R)_{(p)}. \end{array}$$

The map $\mathrm{TC}(\mathbb{Z}_{(p)})_{(p)} \to L_1 \mathrm{TC}(\mathbb{Z}_{(p)})_{(p)} \simeq L_{n+1} \mathrm{TC}(\mathbb{Z}_{(p)})_{(p)}$ is $(-1)$-truncated by the computation of Bökstedt and Madsen [BM94], and Rognes at the prime 2 [Rog99]. It follows that the same holds for $K$-theory. Therefore, Corollary 2.4 applies.

Moreover, Example 2.6 also implies that $\mathrm{gr}^\star_{\mathrm{mot/MU}} \mathrm{TC}(\mathrm{BP}\langle n \rangle)/(p,\dots,v_{n+1})$ has a horizontal vanishing line, so we can apply Proposition 2.5. $\qquad\square$

**Example 2.8.** By [AHW24, Theorem 3.5.7], the mod-$(p,\ldots,v_{n+1})$ syntomic cohomology of (any $\mathbb{E}_1$-algebra form of) the connective Morava $K$-theories $k(n)$ vanishes in degrees $* > m_{n+1} - 1$. As above, we conclude that the localization maps

$$K(k(n))_{(p)} \to L_{n+1}^f K(k(n))_{(p)}$$
$$K(k(n))_{(p)} \to L_{n+1} K(k(n))_{(p)}$$

are $(-1)$-truncated. We will revisit this example in Section 3 without relying on the full computations of [AHW24].

**Example 2.9.** It is shown in [AAR23, Theorem 5.12] that the mod-$(2, \eta, v_1, v_2)$ syntomic cohomology of ko vanishes in degrees $* > 12 = m_{n+1} + 1$. By [BEM17], the mod-$(2, \eta, v_1, v_2^{32})$ syntomic cohomology is the $E_2$-page of a spectral sequence converging to $F_* \mathrm{TC}(\mathrm{ko})$ where $F$ is the quotient of $A(1)$ by a $v_2^{32}$-self-map. The argument of Theorem 2.3 implies that this vanishes in degrees $* > \dim F - 1$, so the localization maps

$$K(\mathrm{ko})_{(2)} \to L_{n+1}^f K(\mathrm{ko})_{(2)}$$
$$K(\mathrm{ko})_{(2)} \to L_{n+1} K(\mathrm{ko})_{(2)}$$

is $(-1)$-truncated.

## 3. Redshift and the Segal Conjecture

To prove truncatedness of $F_* \mathrm{TC}(R)$, for $F$ a type-$(n+2)$ finite complex, [HW22] propose the general strategy of showing separately:

– [*(Weak) Canonical Vanishing*] The canonical map

$$F_* \mathrm{TC}^-(R) \to F_* \mathrm{TP}(R)$$

is zero in large degrees.

– [*Segal Conjecture*] The frobenius

$$F_* \mathrm{THH}(R) \to F_* \mathrm{THH}(R)^{\mathrm{t}C_p}$$

is an isomorphism in large degrees (and thus the same for $F_* \mathrm{TC}^-(R) \to F_* \mathrm{TP}(R)$).

In [HW22, Theorem 3.3.2(e)], however, it is also shown that under certain extra assumptions, truncatedness of $F_* \mathrm{TC}(R)$ can be reduced to the Segal Conjecture alone. This can be applied, for instance, to $\mathrm{BP}\langle n \rangle$-algebras as soon as one knows redshift for $\mathrm{BP}\langle n \rangle$. In this section, we record analogous effective statements for syntomic cohomology (Theorem 3.6). The requisite feature will be an element $\sigma^2 v_{n+1} \in \mathrm{gr}^\star_{\mathrm{mot}/A} \mathrm{THH}(R)/(p,\ldots,v_n)$ which behaves like a version of the *Bökstedt class* described in [BHLS23, §2.3]. In fact, all that is required of this element is that it be a nonzerodivisor (Proposition 3.5). This tells us that redshift reduces to the Segal conjecture as soon as one checks conditions on *only* the level of THH, rather than the more difficult $\mathrm{TC}^-$ or TP.

As a corollary (Theorem 3.13), we give a short proof of effective redshift for the Morava $K$-theories $k(n)$ without computing $\mathrm{gr}^\star_{\mathrm{mot}/\mathrm{MU}} \mathrm{TC}(k(n))/(p,\ldots,v_{n+1})$. As noted previously, this full computation is already known due to the recent work [AHW24]; nevertheless, we believe that it still may be of conceptual interest that Quillen-Lichtenbaum is "mainly" a consequence of the Segal conjecture.

**Construction 3.1.** Let $R$ be an $\mathbb{E}_3$-algebra over a chromatically $p$-quasisyntomic $\mathbb{E}_\infty$-ring $A$, and suppose that:

(i) the map $\mathrm{MU}_*A \to \mathrm{MU}_*R$ is a *quasiregular quotient* in the sense of [HRW25, Definition 4.1.1], that is, the cotangent complex $\mathrm{L}_{\mathrm{MU}_*R/\mathrm{MU}_*A}$ has Tor-amplitude concentrated in degree 1. (Here we temporarily ignore the grading on $\mathrm{MU}_*R$ and $\mathrm{MU}_*A$.)

(ii) the elements $p, \ldots, v_n$ form a regular sequence in $\mathrm{MU}_*R$.

Recalling the results of [HRW25, §3-4], we can compute $\mathrm{fil}^*_{\mathrm{mot}/A}\mathrm{THH}(R)$, together with the ($p$-completed) canonical and frobenius maps, by constructing a suitable quasiregular quotient $S \to A \otimes \mathrm{MW}$ and then applying descent along $\mathrm{THH}(A) \to \mathrm{THH}(A \otimes \mathrm{MW}/S)$:

$$\mathrm{fil}^*_{\mathrm{mot}/A}\mathrm{THH}(R) \simeq \lim_{\Delta} \tau_{\geq 2*}\mathrm{THH}(R \otimes \mathrm{MW}^{\otimes \bullet + 1}/S^{\otimes \bullet + 1})$$

There is a filtration of $\mathrm{THH}(R \otimes \mathrm{MW}/S)$ with associated graded $\mathrm{THH}(\mathrm{MW}_*R/S_*)$, which in turn is equipped with the HKR filtration whose associated graded is $\mathrm{LSym}\left(\Sigma\mathrm{L}_{\mathrm{MW}_*R/S_*}\right)$. Now, in the fiber sequence

$$\mathrm{MW}_*R \otimes^{\mathbb{L}}_{\mathrm{MW}_*A} \mathrm{L}_{\mathrm{MW}_*A/S_*} \to \mathrm{L}_{\mathrm{MW}_*R/S_*} \to \mathrm{L}_{\mathrm{MW}_*R/\mathrm{MW}_*A}$$

the last term is a flat base change of $\mathrm{L}_{\mathrm{MU}_*R/\mathrm{MU}_*A}$, and we deduce that the map $S_* \to \mathrm{MW}_*R$ is also a quasiregular quotient. In particular, $\Sigma^{-1}\mathrm{L}_{\mathrm{MW}_*R/S_*}$ is a flat $\mathrm{MW}_*R$-module concentrated in even degrees. Therefore, $\mathrm{LSym}\left(\Sigma\mathrm{L}_{\mathrm{MW}_*R/S_*}\right)$ is even and, using (ii), so is the quotient $\mathrm{LSym}\left(\Sigma\mathrm{L}_{\mathrm{MW}_*R/S_*}\right) \otimes^{\mathbb{L}}_{\mathrm{MW}_*R} \mathrm{MW}_*R/(p, \ldots, v_n)$. It follows from the appropriate spectral sequences that $\mathrm{THH}(R \otimes \mathrm{MW}/S)$ and $\mathrm{THH}(R \otimes \mathrm{MW}/S)/(p, \ldots, v_n)$ are even. We conclude, finally, that the motivic $E_2$-page and the mod-$(p, \ldots, v_n)$ motivic $E_2$ page can be realized by an explicit *cobar complex*:

$$\mathrm{gr}^*_{\mathrm{mot}/A}\mathrm{THH}(R) \simeq \lim_{\Delta} \pi_{2*}\mathrm{THH}(R \otimes \mathrm{MW}^{\otimes \bullet + 1}/S^{\otimes \bullet + 1})$$

$$\mathrm{gr}^*_{\mathrm{mot}/A}\mathrm{THH}(R)/(p, \ldots, v_n) \simeq \lim_{\Delta} \pi_{2*}\mathrm{THH}(R \otimes \mathrm{MW}^{\otimes \bullet + 1}/S^{\otimes \bullet + 1})/(p, \ldots, v_n).$$

**Example 3.2.** The prototypical case is $A = \mathrm{MU}$ and $R = \mathrm{BP}\langle n \rangle$, with the $\mathbb{E}_3$-algebra structure being the main result of [HW22, §2].

**Construction 3.3.** Given a unital algebra $R$ over a base $\mathbb{E}_\infty$-ring $S$, the composite

$$S \longrightarrow R \xrightarrow{\mathrm{id} \otimes 1 - 1 \otimes \mathrm{id}} R \otimes_S R$$

is zero, so we get an induced *suspension map* $\sigma : \Sigma(R/S) \to R \otimes_S R$. As detailed in [HW22, §A.2] and [LL23, §3.1], this generalizes to a double suspension map to Hochschild homology

$$\sigma^2 : \Sigma^2(R/S) \to \mathrm{THH}(R/S)$$

when $R$ is an $\mathbb{E}_1$-algebra. In [HW22, §A.4], it is shown that the double suspension is "undone" by the circle action on THH: that is,

$$t\sigma^2 x = x \mod t^2$$

in $\mathrm{TC}^-(R/S)$.

**Construction 3.4.** Let $R$ and $A$ be as in Construction 3.1, and suppose additionally that

(iii) $v_{n+1}$ is zero in $\mathrm{MU}_*R$ modulo $(p, \ldots, v_n)$.

In the mod-$(p, \ldots, v_n)$ cobar complex for $\mathrm{TC}^-$, we have the relation $t\sigma^2 v_{n+1} = v_{n+1} + xt^2$ from Construction 3.3. In an abuse of notation, we will relabel $\sigma^2 v_{n+1} - xt$ by $\sigma^2 v_{n+1}$, so $t\sigma^2 v_{n+1} = v_{n+1}$ holds on the nose. We now argue that $\sigma^2 v_{n+1}$ is a cocycle in the mod-$(p, \ldots, v_n)$ cobar complex for THH. Starting from the the mod-$(p, \ldots, v_n)$ cobar complex for $\mathrm{TC}^-$, we have

$$0 = d(v_{n+1}) = d(t\sigma^2 v_{n+1}) = \eta_R(t)\eta_R(\sigma^2 v_{n+1}) - t\sigma^2 v_{n+1}.$$

Since $\eta_R(t) = t + O(t^p)$, and $t$ is a nonzerodivisor, we may divide by $t$ to get that $d(\sigma^2 v_{n+1}) = 0 \bmod t$. But modding out by $t$ is exactly how we obtain the mod-$(p, \ldots, v_n)$ cobar complex for THH. In conclusion, we may speak of an element $\sigma^2 v_{n+1}$ in $\mathrm{gr}^*_{\mathrm{mot}/A} \mathrm{THH}(R)/(p, \ldots, v_n)$.

**Proposition 3.5.** *Let $R$ and $A$ be as in* Construction 3.4, *i.e. $R$ is an $\mathbb{E}_3$-algebra over a chromatically $p$-quasisyntomic $\mathbb{E}_\infty$-ring $A$ such that:*

(i) $\mathrm{MU}_* A \to \mathrm{MU}_* R$ *is a quasiregular quotient.*
(ii) *the elements $p, \ldots, v_n$ form a regular sequence in $\mathrm{MU}_* R$.*
(iii) $v_{n+1}$ *is zero in $\mathrm{MU}_* R$ modulo $(p, \ldots, v_n)$.*

*Suppose furthermore that $\sigma^2 v_{n+1}$ is a nonzerodivisor in $\mathrm{gr}^*_{\mathrm{mot}/A} \mathrm{THH}(R)/(p, \ldots, v_n)$. Then $\sigma^2 v_{n+1}$ is invertible under the Frobenius in $\mathrm{gr}^*_{\mathrm{mot}/A} \mathrm{THH}(R)^{tC_p}/(p, \ldots, v_n)$.*

*Proof.* We have the following equation in the zeroth term of the mod-$(p, \ldots, v_n)$ cobar complex for TP:

$$\varphi(v_{n+1}) = \varphi(t\sigma^2 v_{n+1}) = \varphi(t)\varphi(\sigma^2 v_{n+1})$$
$$= ([p](t))\,\varphi(\sigma^2 v_{n+1})$$
$$= (v_{n+1}t^{p^{n+1}} + O(t^{2p^{n+1}}))\,\varphi(\sigma^2 v_{n+1})$$

(where we take the $p$-typical coordinate from BP, and elements are labeled via the canonical map). But as argued in [Hah22, Lemma 2.5], $\varphi(v_{n+1})$ is also equal to just $v_{n+1}$, modulo $(p, \ldots, v_n)$.

We compare the leading coefficients of our expressions for $\varphi(v_{n+1})$. First we rewrite $v_{n+1} = t\sigma^2 v_{n+1}$ where $\sigma^2 v_{n+1}$ is nonzero mod $t$. Let $a$ be the leading coefficient of $\varphi(\sigma^2 v_{n+1})$, i.e.

$$\varphi(\sigma^2 v_{n+1}) = at^k + O(t^{k+1}).$$

We wish to show that $a$ is a unit. Now, $a$ satisfies

$$a\sigma^2 v_{n+1}t^{k+p^{n+1}+1} + O(t^{k+p^{n+1}+2}) = \sigma^2 v_{n+1}t$$

and we see that $a$ must satisfy either $a\sigma^2 v_{n+1} = \sigma^2 v_{n+1}$ (if $k = -p^{n+1}$) or $a\sigma^2 v_{n+1} = 0$ (if $k < -p^{n+1}$) in the mod-$(p, \ldots, v_n)$ cobar complex for THH. We know that $\sigma^2 v_{n+1}$ is a nonzerodivisor after taking cohomology, which in this case is just taking cocycles since we are in weight 0. Therefore, if we can show that $a$ is a cocycle, then in fact $a = 1$ by propagating the cobar relation $a\sigma^2 v_{n+1} = \sigma^2 v_{n+1}$ to cohomology.

Indeed, from the formula $\eta_R(t) = t + O(t^p)$, we deduce that $\eta_R(\varphi(\sigma^2 v_{n+1}))$ has leading coefficient $\eta_R(a)$ in the cobar complex for TP. But since

$$\eta_R(\sigma^2 v_{n+1}) = \sigma^2 v_{n+1} + O(t^{p-1})$$

from Construction 3.4, applying the frobenius to both sides gives

$$\eta_R(\varphi(\sigma^2 v_{n+1})) = \varphi(\sigma^2 v_{n+1}) + O(t^{p^{n+2}})$$
$$= at^k + O(t^{k+1}).$$

In particular, $\eta_R(a) = a \bmod t$, as desired.                              □

**Theorem 3.6.** *Let $R$ and $A$ be as in* Construction 3.4. *Suppose that*

i) $\mathrm{gr}^*_{\mathrm{mot}/A} \mathrm{THH}(R)$ *satisfies the Segal Conjecture such that the Frobenius*

$$\varphi : \mathrm{gr}^*_{\mathrm{mot}/A} \mathrm{THH}(R)/(p, \ldots, v_n) \to \mathrm{gr}^*_{\mathrm{mot}/A} \mathrm{THH}(R)^{tC_p}/(p, \ldots, v_n)$$

*is an isomorphism in degrees $* \geq m_n + d$*

*ii) $\sigma^2 v_{n+1}$ is invertible under the Frobenius in $\mathrm{gr}^\star_{\mathrm{mot}/A} \mathrm{THH}(R)^{\mathrm{t}C_p}/(p, \ldots, v_n)$.*

*Then the localization map*

$$\mathrm{TC}(R) \to L^f_{n+1} \mathrm{TC}(R)$$

*is d-truncated.*

**Remark 3.7.** In fact, the criteria of Theorem 3.6 imply that

$$\mathrm{gr}^\star_{\mathrm{mot}/A} \mathrm{THH}(R)^{\mathrm{t}C_p}/(p, \ldots, v_n) \simeq \mathrm{gr}^\star_{\mathrm{mot}/A} \mathrm{THH}(R)/(p, \ldots, v_n)[(\sigma^2 v_{n+1})^{-1}]$$

with the Frobenius given by inverting $\sigma^2 v_{n+1}$. This is because $\mathrm{fib}(\varphi)$ is bounded by the Segal conjecture, so $\sigma^2 v_{n+1}$ must be nilpotent in the fiber, which implies that $\mathrm{fib}(\varphi)[(\sigma^2 v_{n+1})^{-1}]$ vanishes.

*Proof of Theorem 3.6.* It remains to understand the canonical map

$$\mathrm{can} : \mathrm{gr}^\star_{\mathrm{mot}/A} \mathrm{TC}^-(R)/(p, \ldots, v_{n+1}) \to \mathrm{gr}^\star_{\mathrm{mot}/A} \mathrm{TP}(R)/(p, \ldots, v_{n+1}).$$

Following [HRW25, §6.5], the canonical map is computed via a map of spectral sequences

$$(\mathrm{gr}^\star_{\mathrm{mot}/A} \mathrm{THH}(R)/(p, \ldots, v_n))[t]/t\sigma^2 v_{n+1} \to (\mathrm{gr}^\star_{\mathrm{mot}/A} \mathrm{THH}(R)/(p, \ldots, v_n))[t^\pm]/\sigma^2 v_{n+1}$$

given by inverting $t$. By conditions (i) and (ii), all homotopy classes in the domain must be divisible by $\sigma^2 v_{n+1}$ in the degree range $* \geq (m_n + d + |\sigma^2 v_{n+1}|)$, or equivalently in the degree range $* > (m_{n+1} + d)$. It follows that $\mathrm{gr}^\star_{\mathrm{mot}/A} \mathrm{THH}(R)/(p, \ldots, v_n)$ must vanish mod $\sigma^2 v_{n+1}$ in degrees $* > (m_{n+1} + d)$. We therefore see that in this range, the domain of the canonical map on the above spectral sequences is concentrated in nonnegative $t$-adic filtration, while the codomain is concentrated in negative $t$-adic filtration, so the canonical map must be zero.

Thus, we obtain that $(\varphi - \mathrm{can})$ is an isomorphism in degrees $* > (m_{n+1} + d)$, which allows us to apply Theorem 2.3. □

A useful feature of Theorem 3.6 is that invertibility of $\varphi(\sigma^2 v_{n+1})$ in fact only needs to be checked for some initial $R_0$, after which it automatically extends to any $\mathbb{E}_1$-algebras $R$ over $R_0$. And furthermore, Proposition 3.5 implies that we can get away with checking only $\mathrm{gr}^\star_{\mathrm{mot}/\mathrm{THH}}(R)$ rather than $\mathrm{gr}^\star_{\mathrm{mot}/\mathrm{THH}}(R)^{\mathrm{t}C_p}$.

**Corollary 3.8.** *Let $R_0$ and $A$ be as in Construction 3.4 such that $\varphi(\sigma^2 v_{n+1})$ is invertible in $\mathrm{gr}^\star_{\mathrm{mot}/A} \mathrm{THH}(R)^{\mathrm{t}C_p}/(p, \ldots, v_n)$. Let $R$ be an $R_0$-$\mathbb{E}_1$-algebra which satisfies the Segal Conjecture such that the Frobenius*

$$\varphi : \mathrm{gr}^\star_{\mathrm{mot}/A} \mathrm{THH}(R)/(p, \ldots, v_n) \to \mathrm{gr}^\star_{\mathrm{mot}/A} \mathrm{THH}(R)^{\mathrm{t}C_p}/(p, \ldots, v_n)$$

*is an isomorphism in degrees $* \geq m_n + d$. Then the localization map*

$$\mathrm{TC}(R) \to L^f_{n+1} \mathrm{TC}(R)$$

*is d-truncated.*

*Proof.* We have now that $\mathrm{gr}^\star_{\mathrm{mot}/A} \mathrm{THH}(R)$ is a module over $\mathrm{gr}^\star_{\mathrm{mot}/A} \mathrm{THH}(R_0)$, and the proof of Theorem 3.6 applies verbatim. □

As it happens, Theorem 3.6 is not so useful in practice as of yet, because there are no rings for which we know the Segal Conjecture for its motivically-filtered THH without also having computed its syntomic comology outright. However, there is reason to believe that showing the Segal Conjecture is significantly "easier" than computing syntomic cohomology. Indeed, in the non-motivic setting, [HW22] readily deduce the Segal Conjecture for $\mathrm{THH}(\mathrm{BP}\langle n \rangle)$ from the Segal Conjecture for $\mathrm{THH}(\mathbb{F}_p)$, by using the $\mathbb{F}_p$-Adams filtration on $BP\langle n\rangle$ and checking on associated graded. In particular, this does not require knowing anything about the behavior of the Tate or homotopy fixed point spectral sequences for $TP(BP\langle n\rangle)$ or $TC^-(BP\langle n\rangle)$. It does not seem unreasonable to posit that a suitable theory of the even filtration exists for filtered rings, which would allow one to run a motivic analogue of this argument.

Instead we finish here with the example of the connective Morava $K$-theories, whose syntomic cohomology has been recently computed by [AHW24]. We exhibit a version of the above strategy as a shortcut to deduce effective redshift, as soon as one knows the homotopy groups of $gr^{\star}_{mot/MU}THH(k(n); \mathbb{F}_p)$ as well as the lemma [AHW24, Proposition 3.3.4], both reproduced below:

**Proposition 3.9** ([AHW24, Proposition 2.2.10])**.** *The motivic spectral sequence computing* $THH_*(k(n); \mathbb{F}_p)$ *collapses at the $E_2$-page as*

$$\pi_* gr^{\star}_{mot/MU}THH(k(n); \mathbb{F}_p) \cong \mathbb{F}_p[\sigma^2 v_{n+1}] \otimes \mathbb{F}_p[\mu]/\mu^{p^n} \otimes \Lambda(\lambda_{n+1})$$

*as an $\mathbb{E}_0$-$\pi_* gr^{\star}_{mot/MU}THH(BP\langle n\rangle; \mathbb{F}_p)$-algebra, where $\mu^j$ is in bidegree $(2j, 0)$, and $\lambda_{n+1}$ is in bidegree $(2p^{n+1}, 1)$.*

**Proposition 3.10** ([AHW24, Proposition 3.3.4])**.** *The bigraded homotopy groups*

$$\pi_* gr^{\star}_{mot/MU}THH(k(n))^{tC_p}$$

*are concentrated on the 0-line and the 1-line.*

**Remark 3.11.** Note that the above proposition is quite particular to Morava $k$-theories, so we do not mean to present it here as something that is meant to generalize.

**Corollary 3.12.** *The Segal conjecture holds for $gr^{\star}_{mot/MU}THH(k(n))/v_n$, and the Frobenius is given by inverting $\sigma^2 v_{n+1}$.*

*Proof.* Modding out by $v_n$ introduces at most a nonzero $(-1)$-line, but by construction, elements on motivic $E_2$ pages have weight and degree of the same parity, so there are no room for differentials in the motivic spectral sequence for $THH_*(k(n))^{tC_p}/v_n$. The mod-$v_n$ Frobenius map of motivic spectral sequences therefore has both its source and target collapse, so we only need to check the Segal conjecture non-motivically. Filtering $k(n)$ by its postnikov tower, we obtain a spectral sequence

$$THH_*(\mathbb{F}_p[v_n])/v_n \implies THH_*(k(n))/v_n.$$

The Frobenius on $THH_*(\mathbb{F}_p[v_n])/v_n$ is the map

$$\mathbb{F}_p[x] \otimes \Lambda(\sigma v_n) \to \mathbb{F}_p[x^{\pm}] \otimes \Lambda(\sigma v_n)$$

which is an isomorphism in large degrees, so we deduce that the Segal conjecture holds for $THH_*(k(n))/v_n$.

Finally, we deduce the formula for the Frobenius from Proposition 3.5 and Remark 3.7.
$\square$

**Theorem 3.13.** *The localization maps*

$$K(k(n))_{(p)} \to L^f_{n+1} K(k(n))_{(p)}$$
$$K(k(n))_{(p)} \to L_{n+1} K(k(n))_{(p)}$$

*are $(-1)$-truncated.*

*Proof.* The above discussion implies that

$$\mathrm{gr}^*_{\mathrm{mot}/\mathrm{MU}}\mathrm{THH}(k(n))/(p,\ldots,v_n) \cong \mathbb{F}_p[\sigma^2 v_{n+1}] \otimes \Lambda(\lambda_{n+1}) \otimes \mathbb{F}_p[\mu]/\mu^{p^n} \otimes \Lambda(\sigma p,\ldots,\sigma v_{n-1})$$

and the Frobenius is given by inverting $\sigma^2 v_{n+1}$. Therefore, the Frobenius is an isomorphism on $\mathrm{gr}^*_{\mathrm{mot}/\mathrm{MU}}\mathrm{THH}(k(n))/(p,\ldots,v_n)$ in degrees

$$* > |(\sigma^2 v_{n+1})^{-1} \cdot \lambda_{n+1} \cdot \mu^{p^n-1} \cdot \sigma p \cdot \cdots \cdot \sigma v_n| = m_n - 2.$$

By Theorem 3.6, the localization map

$$\mathrm{TC}(k(n)) \to L^f_{n+1}\mathrm{TC}(k(n))$$

is $(-1)$-truncated. Finally, the result for $K$-theory follows from Corollary 2.4 and Quillen's computation of $K(\mathbb{F}_p)$ [Qui72]. □